\documentclass{birkjour}
\newtheorem{theorem}{Theorem}
 \newtheorem{wm}{WM}

 \newtheorem{jax}{J}
 \numberwithin{equation}{section}

\newtheorem{ax}{A}

\begin{document}

%
%
%
%
%
%
%
%
%

\title[An axiomatic look at a windmill]
 {An axiomatic look at a windmill}

\author{Victor Pambuccian}

\address{%
Division of Mathematical and Natural Sciences (MC 2352)\\
Arizona State University - West Campus\\
P. O. Box 37100\\
Phoenix, AZ 85069-7100\\
U.S.A.}

\email{pamb@asu.edu}







\begin{abstract}
We present the problem stated in intuitive language as problem 2 at
the 52nd International Mathematical Olympiad as a formal statement,
and prove that it is valid in ordered regular incidence planes, the
weakest ordered geometry whose models can be embedded in projective
ordered planes.
\end{abstract}
\maketitle

\section{Introduction}

Proposed by  Geoffrey Smith of the University of Bath, the second
problem on the first day of the 52nd IMO, held in Amsterdam, reads
as follows (see \cite{imo} for the statements and proofs of all problems):\\

Let ${\mathcal S}$ be a finite set of at least two points in the
plane. Assume that no three points of ${\mathcal S}$ are collinear.
A {\em windmill} is a process that starts with a line $l$ going
through a single point $P\in {\mathcal S}$. The line rotates
clockwise about the {\em pivot} $P$ until the first time that the
line meets some other point belonging to ${\mathcal S}$. This point,
$Q$, takes over as the new pivot, and the line now rotates clockwise
about $Q$, until it next meets a point of $S$. This process
continues indefinitely, with the pivot always being a point from
${\mathcal S}$.\\
Show that we can choose a point $P$ in  and a line $l$ going through
$P$ such that the resulting windmill uses each point of ${\mathcal
S}$ as a pivot infinitely many times.\\

As stated, the problem appears on a first reading to be describing a
process in Euclidean geometry or at any rate a geometry with a
metric, as it appears to require the existence of rotations, without
really belonging to Euclidean geometry proper, given its strong
combinatorial flavor.

The aim of this note is to state it as a theorem of ordered regular
incidence planes.

We will proceed as follows: First, we will provide an axiom system
for planar ordered domains (without the lower-dimension axiom), then
we will state the windmill problem inside that formalism, and
provide a proof for it. Next we will introduce axiomatically ordered
regular incidence planes, state the windmill problem in that
formalism, in which it turns out to be a universal statement (i.e.
it does not contain any existential quantifier), and finally provide
the rationale for the validity of the windmill problem inside that
axiom system.

\section{The axiomatic framework for planar ordered domains}

The axiomatic framework is that of a very  general two-dimensional
theory of betweenness, the models of which will be referred to as
{\em planar ordered domains} (see also \cite{pamb}), axiomatized in
terms of {\em points} as individual variables and the strict
betweenness ternary predicate $Z$, with $Z(abc)$ to be read as `$b$
lies between $a$ and $c$' (and the order is strict, i.\ e.\ $b$ is
different from $a$ and $c$), the axiom system consisting of the
axioms A\ref{a1}-A\ref{a5} axiomatizing $\mathcal{L}$, the universal
theory of linear order  (we omit throughout the paper universal
quantifiers in universal sentences):
\begin{ax}\label{a1}
$Z(abc)\rightarrow Z(cba),$
\end{ax}
\begin{ax} \label{a2}
$Z(abc)\rightarrow \neg Z(acb),$
\end{ax}
\begin{ax}\label{a3}
$Z(acb)\wedge Z(abd)\rightarrow Z(cbd),$
\end{ax}
\begin{ax} \label{a4}
$Z(cab)\wedge Z(abd) \rightarrow Z(cbd),$
\end{ax}
\begin{ax} \label{a5}
$c\neq d \wedge Z(abc)\wedge Z(abd) \rightarrow (Z(bcd)\vee
Z(bdc)),$
\end{ax}

the lower-dimension axiom, stating that there are three
non-collinear points, and the Pasch axiom (here $L$ stands for the
collinearity predicate, defined by $L(xyz):\Leftrightarrow
Z(xyz)\vee Z(yzx)\vee Z(zxy)\vee x=y \vee y=z \vee z=x$; although we
do not have the concept of a `line' in our language, we will refer
to lines, saying that the point $p$ lies on the line $\langle a,
b\rangle$, for two distinct points $a$ and $b$, if $p=a$ or $p=b$ or
$Z(apb)$ or $Z(pba)$ or $Z(bap)$):

\begin{ax} \label{pasch}
$(\forall abcde)(\exists f)\, [\neg L(abc)\wedge Z(adc)\wedge \neg
L(ace)\wedge \neg L(edb)$\\
\hspace*{10mm}$\rightarrow (Z(afb)\vee Z(bfc))\wedge L(edf)].$
\end{ax}

Notice that we do not ask the order to be dense or unending, and we
also do not need the lower-dimension axiom, stating the existence of
three non-collinear points.

\section{The windmill statement}

Let $r_m(k)$ denote the remainder of the division of $k$ by $m$, and
let $^{\epsilon}\varphi$ stand for $\varphi$ if $\epsilon=0$, and
for $\neg\varphi$ if $\epsilon=1$. Let $_{\epsilon}(\varphi)$ stand
for $\perp$ (the absurdity sign) if $\epsilon=0$, and for $\varphi$
if $\epsilon=1$. To help with the readability of the formal
statement of the windmill problem, we introduce the following
defined predicates (where $j,k,l\in \{0,1\}$):

\begin{eqnarray*}
\lambda(xyz) & \Leftrightarrow & Z(xyz)\vee Z(yzx)\vee Z(zxy)\\
\delta(abuv) & \Leftrightarrow  & (\exists t)\, \lambda(abt)\wedge
Z(utv)\\
 \pi_{k,l}^j(a,b,p,c) & \Leftrightarrow  &
 (^{r_2(j+k+1)}\delta(acbp)\vee _{r_2(j+l)}(p=c))\wedge
^{r_2(j+l)}\delta(abcp)\\
 & & \wedge \bigwedge_{\{i| a_i\not\in \{a,b,c,p\}\}}\\
 && [(^{(1-k)}\delta(abca_i)\wedge
 ^l\delta(acba_i))\vee (^k\delta(abca_i)\wedge
^{(1-l)}\delta(acba_i))]\\
\end{eqnarray*}

 \begin{figure}
\input epsf
{\hfill\epsfysize 6cm \epsffile{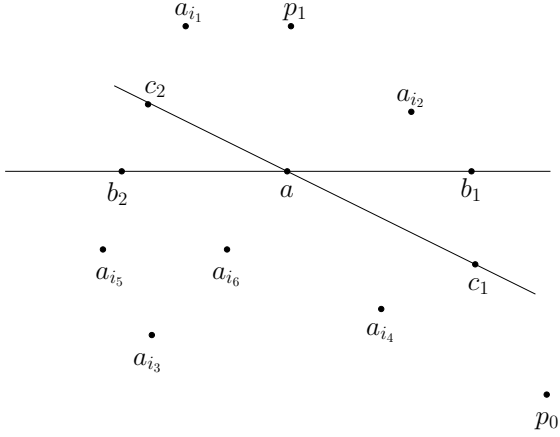} \hfill}
\caption[]{{\small $\pi_{1,1}^i(a,b_1,p_i,c_1)$,
$\pi_{1,0}^i(a,b_1,p_i,c_2)$,  $\pi_{0,1}^i(a,b_2,p_i,c_1)$,
$\pi_{0,0}^i(a,b_2,p_i,c_2)$}}
\end{figure}

Let $\alpha(n)=n(n-1)+1$. Let $K_n=\{f\, |\, f:\{1,2,\ldots,
\alpha(n)\}\rightarrow \{1,2.\ldots, n\},$ $(\exists k(n))\,
k(n)\leq \alpha(n), \mbox{ such that the}$  $\mbox{restriction}$
$\mbox{of } f \mbox{ to } \{1,2,\ldots, k(n)\} \mbox{ is onto,}$
$(\forall i)\, 3\leq i\leq k(n), f(i)\neq f(i-1), f(i)\neq f(i-2),
f(k(n)-1)=f(2)$, $f(k(n))=f(1)\}$, let $^A2$ stand for the set of
all functions from $A$ to $\{0,1\}$, and in case $A=\{1,2,\ldots,
m\}$, let us write $^m2$ for $^A2$. For  $f\in K_n$, we define
$f(0)=r_{n+1}(f(1)+f(2))$.

With these definitions, we are ready to express the windmill theorem
as the following statement:

\begin{wm}
$\bigvee_{1\leq i<j<k\leq n} L(a_ia_ja_k)$\\
\hspace*{11mm}$\vee \bigvee_{f\in K_n, g\in ^{k(n)}2}
\bigwedge_{i=3}^{k(n)}
\pi_{1-g(i-1),g(i)}^{1-g(i-2)}(a_{f(i-1)},a_{f(i-2)},a_{f(i-3)},a_{f(i)}).$
\end{wm}

To see that this is actually the windmill statement of Geoffrey
Smith's problem, notice that $\pi_{k,l}^j(a,b,p,c)$ means that $c$
is the first point in the set $\{c, a_1, \ldots, a_n\}\setminus \{a,
b\}$ the line determined by $a$ and $b$ meets, when it ``rotates"
around $a$ ``clockwise."  The point $p$ is there to fix the sense of
the ``rotation," and the sub- and superscripts $j, k$, and $l$ are
there to tell us  whether $p$ is in the half-plane toward which the
ray $\stackrel{\rightarrow}{ab}$ moves and whether $c$ will be met
by $\stackrel{\rightarrow}{ab}$ or by the opposite ray during its
``rotation" (see Fig.\ 1 for details). If we look at the ``windmill"
process of the original statement by Geoffrey Smith, we notice that
the ``windmill" consists of lines through single points of the given
set $S=\{a_1, \ldots, a_n\}$, and that pivots are obtained during
what one could call ``windmill stops." Instead of focusing on the
``windmill" as composed by lines through single points of $S$, we
decided to look at the ``windmill" as a collection of ``windmill
stops," i.\ e.\ a collection of lines passing through exactly two
points of $S$. The statement ${\bf WM}$ makes is that, given a set
$S$ of points  $a_1, \ldots a_2$, that are such that no three are
collinear, one can find a map $f$ in $K_n$, such that the first
windmill stop goes through  $a_{f(1)}$ and $a_{f(2)}$ ($a_{f(0)}$
being chosen as a point different from $a_{f(1)}$ and $a_{f(2)}$,
playing the role of $p$ in the formula for $\pi$, i.\ e.\
determining the ``clockwise" direction for the entire windmill
process), with windmill pivots $a_{f(i)}$, with $i\in \{1, \ldots,
k(n)\}$, which, given the definition of $K_n$, exhaust $S$, and such
that the last windmill stop has pivot  $a_{f(1)}$, and goes through
$a_{f(1)}$ and $a_{f(2)}$, just like the stop we started with.
During the whole process, the point $p$ determining the orientation
changes (whenever $c$ becomes pivot, and the windmill stop is the
line $ac$, the new point $p$ is chosen to be the point $b$ from the
windmill stop $ab$ that ``rotated" to stop at $c$, i.\ e.\ the $b$
for which $\pi_{k,l}^j(a,b,p,c)$ holds), but the orientation itself,
and thus the ``clockwise" sense of the ``rotation", stays the same.
The reason we cannot stay with the same $p$ is that the windmill can
stop at $p$, i.\ e.\ we can have $c=p$ in $\pi_{k,l}^j(a,b,p,c)$.

\section{The proof}

Since the proof in \cite{imo} is carried out inside the Euclidean
plane, we need to provide a variant thereof that would hold inside
planar ordered domains. We treat the case in which $n$, the number
of points in the set $S$, is even and the case in which it is odd
separately.

If $n$ is even, then let $a_s$ be any vertex of the convex hull of
$S$ (these elementary notions of convex geometry can all be defined
and have the usual properties, as shown in \cite{cop}). Among the
lines $\langle a_s, a_i\rangle$ with $i\in \{1,\ldots, n\}\setminus
\{s\}$ there must be one having an equal number of points on each of
its sides. We denote that particular value of $i$ by $t$, and start
the windmill process with the first pivot in $a_s$ and windmill stop
at $a_sa_t$, i.\ e.\ $f(1)=s$, $f(2)=t$. The first value of $p$ will
be $a_{r_{n+1}(f(1)+f(2))}$ (choosing the index to be
$r_{n+1}(f(1)+f(2))$ had no other function than making sure it is
different from $f(1)$ and $f(2)$), and we'll choose $g(1)=1$,
$g(2)=0$. We define the direction
$\stackrel{\longrightarrow}{a_{f(1)}a_{f(2)}}$ as the {\em East}
(and thus $\stackrel{\longrightarrow}{a_{f(2)}a_{f(2)}}$ as the {\em
West}, the half-plane determined by $\langle a_{f(1)},
a_{f(2)}\rangle$ in which $a_{f(0)}$ lies as the {\em Southern}
half-plane. These directions change in the course of the windmill
process as follows: if  $\pi_{k,l}^j(a,b,p,c)$, then
$\stackrel{\longrightarrow}{ab}$ and
$\stackrel{\longrightarrow}{ac}$ point in the same ``direction" (i.\
e.\ both East or both West) if $k=l$ and in different directions if
$k\neq l$, whereas the half-plane determined by $ab$ in which $p$
lies has the same name as the half-plane determined by $ac$ in which
$b$ lies if $k=j$, and the opposite name if $k\neq j$. The Southern
half-plane determined by a windmill stop $ab$ will be denoted by
$\sigma_{ab}$, and the Northern half-plane by $\nu_{ab}$

We now look at the possible changes in the difference
$\delta(ab)=N(ab)-S(ab)$ between the number $N(ab)$ of points in
$\nu_{ab}$ and the number $S(ab)$ of points in  $\sigma_{ab}$ during
the windmill process. We want to show that $\delta$ can take only
the values $0$ and $2$.

The next stop will be a point $a_{f(3)}$, with the property that
there is no point in $S$ between the rays
$\stackrel{\longrightarrow}{a_{f(1)}a_{f(2)}}$ and
$\stackrel{\longrightarrow}{a_{f(1)}a_{f(3)}}$. Given that
$a_{f(1)}$ is a vertex of the convex hull of $S$, and that the sense
of the ``rotation about $a_{f(1)}$"  is ``clock-wise", i.\ e.\
towards the Southern half-plane, $a_{f(3)}$ must lie in
$\sigma_{a_{f(1)}a_{f(2)}}$, thus $\delta(a_{f(1)}a_{f(3)})$ must be
$2$, as the Northern half-plane gains a point, namely $a_{f(2)}$,
and the Southern half-plane loses one, namely $a_{f(3)}$.

Point $a_{f(3)}$ becomes a pivot during the next stage of the
windmill process, and at the next stop, i.\ e.\ when line $\langle
a_{f(1)}, a_{f(3)}\rangle $ ``rotates about $a_{f(3)}$" into
$\langle a_{f(3)}, a_{f(4)}\rangle $ (where $a_{f(4)}$ is the point
in $S$ which lies either in $\nu(a_{f(3)}a_{f(1)})$ and for which
there is no point in $S$ between
$\stackrel{\longrightarrow}{a_{f(3)}a_{f(1)}}$ and
$\stackrel{\longrightarrow}{a_{f(3)}a_{f(4)}}$, nor between the rays
opposite to the above two, or else lies in
$\sigma(a_{f(3)}a_{f(1)})$, and no point in $S$ lies between the ray
opposite to $\stackrel{\longrightarrow}{a_{f(3)}a_{f(4)}}$ and
$\stackrel{\longrightarrow}{a_{f(3)}a_{f(1)}}$, nor between the rays
opposite to the above two), point $a_{f(1)}$ will be in
$\sigma_{a_{f(3)}a_{f(4)}}$, and we distinguish two cases: (i)
$a_{f(4)}$ is in $\sigma_{a_{f(1)}a_{f(3)}}$ and (ii) $a_{f(4)}$ is
in $\nu_{a_{f(1)}a_{f(3)}}$. In case (i),
$\sigma_{a_{f(3)}a_{f(4)}}$ both gains and loses a point when
compared to $\sigma_{a_{f(1)}a_{f(3)}}$, and thus $\delta$ stays the
same, namely $2$, and we are back to the situation we were in at the
windmill stop $a_{f(1)}a_{f(3)}$, with the pivot, $a_{f(4)}$, East
of the other point, $a_{f(3)}$, of the windmill stop, with $\delta$
taking the value $2$.

In case (ii),  $\sigma_{a_{f(3)}a_{f(4)}}$ will have one point,
$a_{f(1)}$, more than $\sigma_{a_{f(1)}a_{f(3)}}$, and
$\nu_{a_{f(3)}a_{f(4)}}$ will have one point, $a_{f(4)}$, less than
$\nu_{a_{f(1)}a_{f(3)}}$, so $\delta$ will become $0$ for the
windmill stop $a_{f(3)}a_{f(4)}$, and we  are back to a
configuration of the type we started with at windmill stop
$a_{f(1)}a_{f(2)}$, with the pivot, $a_{f(4)}$, to the West of the
other point of the windmill stop, $a_{f(3)}$. However, there is one
difference:  $a_{f(4)}$ no longer needs to be a vertex of the convex
hull of $S$, and so we can no longer say, as in the case of the
windmill stop $a_{f(1)}a_{f(2)}$, that the next point the windmill
will meet, while ``rotating about $a_{f(4)}$", will lie in
$\sigma_{a_{f(4)}a_{f(3)}}$. It can lie in either the Northern and
the Southern half-plane, and so a new situation can appear, that in
which $a_{f(5)}$ lies in $\nu_{a_{f(3)}a_{f(4)}}$. There is thus,
one situations that still needs to be analyzed in complete
generality: the pivot $a$ is West of $b$ at the windmill stop $ab$,
with $\delta(ab)=0$. The case left open is that in which $c$,  the
first point line $\langle a, b\rangle$ meets while ``rotating about
$a$", lies in $\nu_{ab}$. In that case  $\nu_{ac}$ both gains and
loses a point with respect to $\nu_{ab}$ (it gains $b$ and loses
$c$), so we are back into the previous situation, i.\ e.\ the pivot
$c$ lies to the West of the other point, $a$, of the windmill stop
and $\delta(ac)=0$.

We conclude that, during the entire windmill process, we have only
the following two possible situations for the windmill stop $ab$:
either $\delta(ab)=0$, and the pivot $a$ lies to the West of $b$, or
else $\delta(ab)=2$, and the pivot $a$ lies to the East of $b$.

At every stage of the windmill process, the {\em imprint} of the
``East" on the convex hull of $S$, i.\ e.\ the point where the
Eastward pointing ray of that windmill stop intersects the convex
hull of $S$, moves in clockwise direction towards $a_{f(1)}$, except
when the pivot is a vertex of the convex hull of $S$, in which case
the imprint stays put for that one step in the process, but will
have to move on in the next, as the pivot changes at that step.
After a finite number of steps, in fact in no more than $n(n-1)/2$
steps (given that this is the total number of lines that can be
formed ny joining two points in $S$, and thus the upper bound on the
number of windmill stops), the imprint will be for the last time in
$\sigma_{a_{f(1)}a_{f(2)}}$, in the sense that at the next windmill
stop the imprint of the ``East" will have to be either $a_{f(1)}$ or
be in $\nu_{a_{f(1)}a_{f(2)}}$. However, it cannot be in
$\nu_{a_{f(1)}a_{f(2)}}$, for $a_{f(1)}$ would have lied between the
two windmill stops, contradicting the definition of $\pi$.

The pivot $a$ of the last windmill stop $\langle a, b\rangle$ for
which the Eastward imprint is in $\sigma_{a_{f(1)}a_{f(2)}}$ cannot
lie in $\sigma_{a_{f(1)}a_{f(2)}}$, for if it did, then  $\langle a,
a_{f(1)}\rangle$  would be the next windmill stop, and
$\delta(aa_{f(1)})$ would have to be $\leq -1$, contradicting the
fact that $\delta$ takes only non-negative values.

If $a$ lies on $\langle a_{f(1)}, a_{f(2)}\rangle$, then $a$ would
have to be $a_{f(2)}$, and, since the next windmill stop, after
$\langle a_{f(2)}, b\rangle$ is $\langle a_{f(1)}, a_{f(2)}\rangle$,
we would be back in the starting position, with $\langle a_{f(1)},
a_{f(2)}\rangle$ as windmill stop and $a_{f(1)}$ as pivot, and thus
we'd be done.

If $a$ lies in $\nu_{a_{f(1)}a_{f(2)}}$, then, since
$\delta(aa_{f(1)})$ would have to be positive, and since $2$ is the
only positive value it is allowed to take, there can be no point in
$S$ lying between the rays
$\stackrel{\longrightarrow}{a_{f(1)}a_{f(2)}}$ and
$\stackrel{\longrightarrow}{a_{f(1)}a}$, so the next windmill stop,
after $\langle a_{f(1)}, a\rangle$, has to be $\langle a_{f(1)},
a_{f(2)}\rangle$. In that case, we are not quite back where we
started from, for although the windmill stop is $\langle a_{f(1)},
a_{f(2)}\rangle$, the pivot is $a_{f(2)}$, not $a_{f(1)}$ the way it
was at the start. However, by the same reasoning that showed us
that, in case we start with a line which has the same number of
points on each of its sides as first windmill stop, we arrive in at
most $n(n-1)/2$ steps back to itself, we conclude that we'll be
either back to $\langle a_{f(1)}, a_{f(2)}\rangle$ with $a_{f(2)}$
as pivot, in which case starting with  $\langle a_{f(1)},
a_{f(2)}\rangle$ as windmill stop and with $a_{f(2)}$ as pivot, we
arrive back to the same windmill stop and pivot in at most
$n(n-1)/2$ steps, or else we'll be back to $\langle a_{f(1)},
a_{f(2)}\rangle$ with $a_{f(1)}$ as pivot, in which case starting
with  $\langle a_{f(1)}, a_{f(2)}\rangle$ as windmill stop and with
$a_{f(1)}$ as pivot, we arrive back to the same windmill stop and
pivot in at most $n(n-1)$ steps. That each point in $S$ must have
become a pivot by the time the windmill stop returns for the first
time to $\langle a_{f(1)}, a_{f(2)}\rangle$ is easily seen by
noticing that, at this stage $\nu(a_{f(1)}a_{f(2)})$ has become what
used to be $\sigma(a_{f(1)}a_{f(2)})$ at the start of the process,
and that a point in $S$ can move from the Southern to the Northern
half-plane only by having been touched by a windmill stop.

In case $n$ is odd, say $n=2k+1$, we let $a_s$ be any vertex of the
 convex hull of $S$ and choose among the lines $\langle a_s,
a_i\rangle$ with $i\in \{1,\ldots, n\}\setminus \{s\}$ the one
having $k+1$ points on one side and $k$ points on the other side.
Just like in the case in which $n$ is even, we denote that
particular value of $i$ by $t$, and start the windmill process with
the first pivot in $a_s$ and windmill stop at $\langle a_s,
a_t\rangle$, i.\ e.\ $f(1)=s$, $f(2)=t$. We now distinguish two
possibilities:  the first value of $p$ chosen to start the windmill
process, i.\ e.\ $a_{r_{n+1}(f(1)+f(2))}$, can be (i) in the
half-plane with $k$ points or (ii) in the half-plane with $k+1$
points. In case (i), we notice, as in the $n$ even case treated
earlier, that during the windmill process $\delta$ can take on only
two values, namely $1$ and $3$. If we follow the path of the imprint
of the West on the convex hull of $S$, we notice that it moves in
the ``clockwise" direction (i.\ e.\ moving inside
$\nu_{a_{f(1)}a_{f(2)}}$) at every step of the process, unless the
pivot is a vertex of the convex hull of $S$, in which case it rests
for one step of the windmill process, but will have to move
afterwards. After a finite number of steps, no more than twice the
total number of lines that can be formed by two points in $S$ (since
each such line has two directions that can become the ``West"
direction during the windmill process), the imprint of the West has
to come back to its original location, $a_{f(1)}$ (it cannot jump
over it, the reason being the same as in the $n$ already discussed
even case). Now, the other point $a$ of the windmill stop  $\langle
a, a_{f(1)}\rangle$  we arrive at, when the imprint of the West is
back at $a_{f(1)}$ must be $a_{f(2)}$. To see this, notice that, if
$a$ were in $\sigma_{a_{f(1)}a_{f(2)}}$, then $\delta(aa_{f(1)})$
would have to be $\geq 3$, so at the next windmill stop  $\langle
a_{f(1)}, b\rangle$, we'd have $\delta(a_{f(1)}b)\geq 5$, which is
not possible. If $a$ were in $\nu_{a_{f(1)}a_{f(2)}}$, then
$\delta(aa_{f(1)})$ would have to be $\leq -1$, which is impossible,
and thus $a=a_{f(2)}$. Case (ii) is treated analogously, by noticing
that throughout the windmill process $\delta$ takes on only the
values $-1$ and $1$. For reasons similar to those mentioned in the
$n$ even case, each point in $S$ must have become a pivot during the
windmill process.

\section{Validity in ordered regular incidence planes}

There is a weaker axiom system, for {\em ordered regular incidence
planes} from which ${\bf WM}$ can be derived. It cannot be expressed
in terms of {\em points} and $Z$, as it is based on the notion of
{\em sides} of a line in a plane, put forward by Sperner in
\cite{spe49a}, from which $Z$ can be defined, but which cannot, in
general, be defined in terms of $Z$. It can be expressed in a
two-sorted language, with variables for {\em points} (to be
represented by lower-case Latin characters) and for {\em lines} (to
be represented by lower-case Gothic characters), with two relation
symbols, $I$, with $I(a{\mathfrak g})$ to be read as `point $a$ is
incident with line ${\mathfrak g}$', and $D$, with $D(a{\mathfrak
g}b)$ to be read as `the points $a$ and $b$ lie on different sides
of line ${\mathfrak g}$'.  With $\delta(ab{\mathfrak g}{\mathfrak
h}):\Leftrightarrow [(D(a{\mathfrak g}b)\wedge D(a{\mathfrak
h}b))\vee (\neg D(a{\mathfrak g}b)\wedge \neg D(a{\mathfrak h}b))]$
and $^{\epsilon}\delta$ standing for $\delta$ if $\epsilon=1$ and
for $\neg \delta$ if $\epsilon=0$, the axioms are (see
\cite{jouss}):

\begin{jax} \label{sp1}
$(\forall ab)(\exists^{=1} {\mathfrak g})\, a\neq b \rightarrow
I(a{\mathfrak g})\wedge I(b{\mathfrak g}),$
\end{jax}

\begin{jax} \label{sp1'}
$(\forall {\mathfrak g})(\exists a_1a_2a_3a_4 )\, \bigwedge_{1\leq
i<j\leq 4} a_i\neq a_j \wedge \bigwedge_{i=1}^4 I(a_i{\mathfrak g})$
\end{jax}

\begin{jax} \label{sp2}
$(\exists abc)(\forall {\mathfrak g})\, \neg(I(a{\mathfrak g})\wedge
I(b{\mathfrak g})\wedge I(c{\mathfrak g})),$
\end{jax}

\begin{jax} \label{sp3}
$D(a{\mathfrak g}b)\rightarrow \neg I(a{\mathfrak g}),$
\end{jax}

\begin{jax} \label{sp4}
$D(a{\mathfrak g}b)\rightarrow D(b{\mathfrak g}a),$
\end{jax}

\begin{jax} \label{pa}
$\neg I(c{\mathfrak g})\wedge D(a{\mathfrak g}b)\rightarrow
(D(a{\mathfrak g}c)\vee D(b{\mathfrak g}c)),$
\end{jax}

\begin{jax} \label{sp6}
$\neg (D(a{\mathfrak g}b)\wedge D(b{\mathfrak g}c) \wedge
D(c{\mathfrak g}a)),$
\end{jax}

\begin{jax} \label{z=1}
$[\bigwedge_{1\leq i<j\leq 4}a_{i}\neq a_{j}\wedge {\mathfrak
h}_i\neq {\mathfrak h}_{j}\wedge \bigwedge_{i=1}^4 I(a_i{\mathfrak
h}_i)\wedge {\mathfrak h}_i\neq {\mathfrak g}\wedge
$\\
\hspace*{14mm}$((\bigwedge_{i=1}^4 I(a_i{\mathfrak g}))\vee
(\bigwedge_{i=1}^4 I(o{\mathfrak
h}_i)))]$\\
\hspace*{14mm}$\rightarrow[\bigvee_{\stackrel{\epsilon_i\in
\{0,1\}}{\epsilon_1+\epsilon_2+\epsilon_3}=2}$
$^{\epsilon_1}\delta(a_3a_4{\mathfrak h}_1{\mathfrak h}_2)\wedge$
$^{\epsilon_2}\delta(a_2a_4{\mathfrak h}_1{\mathfrak h}_3)\wedge$
$^{\epsilon_3}\delta(a_2a_3{\mathfrak h}_1{\mathfrak h}_4)].$
\end{jax}

J\ref{pa} is a weak variant of Pasch's axiom, stating that if a line
${\mathfrak g}$ does not pass through any of the points  $a, b,$ and
$c$, and $a$ and $b$ are on different sides of ${\mathfrak g}$ then
so are at least one of the pairs $\{a, c\}$ and $\{b, c\}$.
J\ref{sp6} is a variant of Pasch's theorem, stating that a line
cannot separate all three pairs $\{a,b\}$,  $\{b,c\}$, and
$\{c,a\}$. One of its special cases, when $a=b=c$, implies that $a$
and $b$ can be on different sides of   ${\mathfrak g}$ only if
$a\neq b$. That these  versions are called ``weak" stems from the
fact that, if a line ${\mathfrak g}$ separates the points $a$ and
$b$, it no longer means that there is a point on ${\mathfrak g}$
which is between $a$ and $b$. Indeed, the line ${\mathfrak g}$ and
the line determined by $a$ and $b$ may have no point in common (a
simple example is provided by the submodel of the ordered affine
plane over ${\mathbb Q}$ whose points have coordinates whose
denominators are powers of $2$, with the plane separation relation
inherited from the ordered affine plane over ${\mathbb Q}$; see
\cite{pambfer} for other examples). The meaning of J\ref{z=1} is
best understood in terms of the notion of {\em separation} $//$
(with $ab//cd$ to be read as the point-pair $(a,b)$ separates the
point-pair $(c.d)$), defined by

\begin{eqnarray} \label{defsep}
a_1a_2//a_3a_4 & :\Leftrightarrow & (\exists {\mathfrak g}{\mathfrak
h}{\mathfrak k})\, \bigwedge_{i=1}^4 I(a_i{\mathfrak g})\wedge
\bigwedge_{1\leq i<j\leq 4} a_i\neq a_j\wedge I(a_1{\mathfrak
h}\wedge I(a_2{\mathfrak k})\nonumber\\
& & \wedge {\mathfrak h}\neq {\mathfrak g}\wedge {\mathfrak k}\neq
{\mathfrak g}\wedge  \neg \delta(a_3a_4{\mathfrak h}{\mathfrak k}).
\end{eqnarray}

One part of it (corresponding to the $\bigwedge_{i=1}^4
I(a_i{\mathfrak g})$ disjunct) states that, if $a_1, a_2, a_3, a_4$
are four different collinear points, then exactly one of the
separation relations $a_1a_2//a_3a_4$, $a_1a_3//a_2a_4$,
$a_1a_4//a_2a_3$ holds. Its other part (corresponding to the
$\bigwedge_{i=1}^4 I(o{\mathfrak h}_i)$ disjunct) is the dual
statement (in the sense of projective geometry).

Joussen \cite{jouss} showed that any model ${\mathfrak M}$ of
J\ref{sp1}-J\ref{z=1} can be embedded in a projective ordered plane
${\mathfrak P}$, whose separation relation $//_{\mathfrak P}$ is an
extension of the separation relation $//_{\mathfrak M}$, defined in
${\mathfrak M}$ terms of $I_{\mathfrak M}$ and $D_{\mathfrak M}$ by
(\ref{defsep}).

The windmill statement ${\bf WM}$ remains the same, if we change the
definition of the defined notions $L$ and $\delta$ occurring in it
to:

\begin{eqnarray*}
\delta(abuv) & \Leftrightarrow & (\exists  {\mathfrak g})\, a\neq b
\wedge I(a{\mathfrak g})\wedge I(b{\mathfrak g})\wedge D(u{\mathfrak
g}v),\\
L(abc) & \Leftrightarrow & (\exists  {\mathfrak g})\, (I(a{\mathfrak
g})\wedge I(b{\mathfrak g}\wedge I(c{\mathfrak g}))\vee a=b\vee b=c
\vee c=a.
\end{eqnarray*}

To see that ${\bf WM}$ is true in ordered regular incidence planes,
suppose ${\bf WM}$ were not derivable from J\ref{sp1}-J\ref{z=1}.
Then there would have to exist a model ${\mathfrak M}$ of
J\ref{sp1}-J\ref{z=1} in which ${\bf WM}$ is false, i.\ e.\ in which
$\neg{\bf WM}$ holds. Now notice that $\neg{\bf WM}$  can be
expressed as an existential statement in the following way:

\begin{eqnarray*}
& & (\exists a_i)_{1\leq i \leq n}(\exists {\mathfrak
g}_{ij})_{1\leq i<j\leq n}\, \bigwedge_{1\leq i<j\leq n} a_i\neq
a_j\wedge \bigwedge_{1\leq i<j<k\leq n}
I(a_i{\mathfrak g}_{ij})\wedge I(a_j{\mathfrak g}_{ij})\\
& & \wedge \bigwedge_{(i,j)\neq (k,l), 1\leq i<j\leq n, 1\leq
k<l\leq n} {\mathfrak g}_{ij}\neq {\mathfrak g}_{kl}\\
& & \wedge \bigwedge_{f\in K_n, g\in ^{k(n)}2}
\left[\bigvee_{i=3}^{k(n)}
\neg\pi_{1-g(i-1),g(i)}^{1-g(i-2)}(a_{f(i-1)},a_{f(i-2)},a_{f(i-3)},a_{f(i)})\right].
\end{eqnarray*}

in which the $\delta(a_ia_juv)$ occurring in the $\pi$'s are just
$D(u{\mathfrak g}_{ij}v)$.

Then $\neg{\bf WM}$, as an existential statement, would have to hold
in the ordered projective plane ${\mathfrak P}$, in which
${\mathfrak M}$ can be embedded, as well. If we remove from the
projective plane ${\mathfrak P}$ a line which does not contain any
of the points $a_i$ which $\neg{\bf WM}$ claims to exist in
${\mathfrak M}$, such that the windmill process does not close
regardless of the choice of its starting position, we obtain a model
${\mathfrak N}$ of A\ref{a1}-A\ref{pasch} in which $\neg{\bf WM}$
holds, a contradiction.

\begin{theorem} \label{joussen}
$\{$J\ref{sp1}-J\ref{sp1'}, J\ref{sp3}-J\ref{z=1}$\}$  $\vdash {\bf
WM}$, where ${\bf WM}$ is expressed in terms of points, lines, $I$,
and $D$.
\end{theorem}

By defining $Z$ in terms of $I$ and $D$ by

\begin{equation} \label{DEFZ}
Z(abc):\Leftrightarrow (\exists {\mathfrak g}{\mathfrak h})\,
{\mathfrak h}\neq {\mathfrak g}\wedge I(a{\mathfrak g})\wedge
I(b{\mathfrak g})\wedge I(c{\mathfrak g})\wedge I(b{\mathfrak
h})\wedge D(a{\mathfrak h}c),
\end{equation}

one can compare the set of $Z$-consequences of the axiom system
$\{$J\ref{sp1}-J\ref{sp1'}, J\ref{sp3}-J\ref{z=1}$\}$  to
$\{$A\ref{a1}-A\ref{pasch}$\}$.  It turns out that the $Z$ defined
by (\ref{DEFZ}) satisfies A\ref{a1}-A\ref{a5}, but does not need to
satisfy A\ref{pasch}. On the other hand, J\ref{sp1'} does not follow
from $\{$A\ref{a1}-A\ref{pasch}$\}$. Although formally incomparable,
intuitively $\{$J\ref{sp1}-J\ref{sp1'}, J\ref{sp3}-J\ref{z=1}$\}$ is
the weaker axiom system.

\end{document}